\documentclass[12pt]{article}

\usepackage{color}
\usepackage{dsfont}

\newtheorem{theorem}{Theorem}

\newcommand{\mymod}[3]{$#1 \equiv #2 \,\, (\!\!\!\!\mod #3)$}

\usepackage{epsfig}
\usepackage{url}
\usepackage{amssymb}
\usepackage{amsmath}
\usepackage{amsfonts}

\def\bZ{{\bf Z}}

\def\Dbar{\leavevmode\lower.6ex\hbox to 0pt{\hskip-.23ex \accent"16\hss}D}

\begin{document}

{\bf\LARGE
\begin{center}
Symmetric Hadamard matrices of order 116 and 172 exist
\end{center}
}

{\Large
\begin{center}
Olivia Di Matteo\footnote{ University of Waterloo, Department of Physics and Astronomy, Institute for Quantum Computing
e-mail: \url{odimatte@uwaterloo.ca}},
Dragomir {\v{Z}}. {\Dbar}okovi{\'c}\footnote{
University of Waterloo, Department of Pure Mathematics,
Institute for Quantum Computing,
Waterloo, Ontario, N2L 3G1, Canada
e-mail: \url{djokovic@uwaterloo.ca}},
Ilias S. Kotsireas\footnote{Wilfrid Laurier University, Department of Physics \& Computer Science,
Waterloo, Ontario, N2L 3C5, Canada,
e-mail: \url{ikotsire@wlu.ca}}
\end{center}
}

\begin{abstract}
We construct new symmetric Hadamard matrices of orders $92,116$, and $172$. While the existence of those of order $92$ was known since
1978, the orders $116$ and $172$ are new. Our construction is based on a recent new combinatorial array discovered by N. A. Balonin and J. Seberry. For order $116$ we used an adaptation of an algorithm for parallel collision search. The adaptation pertains to the modification of some aspects of the algorithm to make it suitable to solve a 3-way matching problem. We also point out that a new infinite series of symmetric Hadamard matrices arises by plugging into the GP array the matrices constructed by Xia, Xia, Seberry, and Wu in 2005.
\end{abstract}

\section{Introduction}
A Hadamard matrix of order $n$ is an $n \times n$ matrix
$H=(h_{ij})$ with elements $\pm 1$ such that
$H H^T = H^T H = n I_n$, where $I_n$ is the $n \times n$ identity
matrix and $T$ stands for transposition. If in addition $H$ is symmetric, i.e. $h_{ij} = h_{ji}$, then
it is called a symmetric Hadamard matrix. A list of currently open cases for symmetric Hadamard matrices appears
in \cite{CraigenKharaghani:2007}, page 277, Table 1.52, and states that the only $12$ odd values of $v < 100$
for which a symmetric Hadamard matrix of order $4v$ is not known to exist are:
$$
    23, 29, 39, 43, 47, 59, 65, 67, 73, 81, 89, 93.
$$
However, symmetric conference matrices of order $46$ were
constructed by R. Mathon \cite{Mathon:1978} in 1978, and it is well known that this implies the existence of symmetric Hadamard matrices of order $2\cdot46=92$. By modifying Mathon's construction,
Balonin and Seberry \cite{Balonin:Seberry:2014} have constructed two more symmetric conference matrices of order $46$ inequivalent to those of Mathon.

In this paper we construct symmetric Hadamard matrices of orders
$92,116,172$. All of them are constructed by using the GP array of Balonin and Seberry \cite{Balonin:Seberry:2015}:
\begin{equation}
    GP =
    \left(
    \begin{array}{cccc}
    A  & BR    & CR    & DR \\
    CR & D^TR  & -A    & -B^T R \\
    BR & -A    & -D^TR & C^T R \\
    DR & -C^TR & B^TR  & -A \\
    \end{array}
    \right),
    \label{GP-array}
\end{equation}
where $R$ is the back-diagonal matrix obtained from the identity matrix by reversing the order of rows. In order to obtain a
symmetric Hadamard matrix of order $4v$ from this array we need four circulant $\{\pm1\}$-matrices (also known as {\em binary matrices}) $A,B,C,D$ of order $v$ such that

(i) $AA^T+BB^T+CC^T+DD^T=4vI_v$;

(ii) $A^T=A$ and $B=C$.

\noindent Such quadruples $[A,B,C,D]$ can be constructed from suitable
difference families, also known as {\em supplementary difference sets (SDS)}, in the cyclic group $\bZ_v$ consisting of four blocks. The authors of \cite{Balonin:Seberry:2015} refer to
these quadruples as {\em propus matrices} and to the corresponding symmetric Hadamard matrix as a {\em propus-Hadamard matrix}.

In the case $v=23$, our construction is quite different from those in \cite{Mathon:1978} and \cite{Balonin:Seberry:2014}.
In the case $v=29$ we use a method for parallel collision search.  We had to adapt this method in order to be able to apply it to the problem of searching for suitable quadruples $[A,B,C,D]$.
In the case $v=43$ we in fact construct D-optimal matrices
(see e.g. \cite{DK:JCD:2012})
of order $86$ from two binary circulants $A,D$ with $A^T=A$.
The blocks $B$ and $C=B$ are provided by the Paley difference set in $\bZ_{43}$.

The smallest order for which the existence question for symmetric Hadamard matrices is still undecided is now
$4 \cdot 39=156$.

Moreover we give a new infinite series of symmetric Hadamard matrices derived from the series of Hadamard matrices
constructed in \cite[Theorem 3]{XXSW:2005}.

\section{Some infinite series of symmetric Hadamard matrices}

We summarize some results pertaining to the existence of infinite series of symmetric Hadamard matrices. The summary is far from being exhaustive.

The following result is proved in
\cite[Corollary 4.6.5]{IoninShrikhande:CUP:2006}.

\begin{theorem}
If \mymod{q}{3}{4} is a prime power and $q+2$ is a prime power, then there exists a symmetric conference matrix of order
$q^2(q+2)+1$ and a symmetric Hadamard matrix of order
$2q^2(q+2)+2$.
\end{theorem}

Note that for $q=3$ one obtains a symmetric Hadamard matrix of order $92$.

\noindent The following result is mentioned in
\cite[Theorem 1.48, p. 277] {CraigenKharaghani:2007}.
\begin{theorem}
If $n+1$ and $n-1$ are both odd prime powers, then there exists a
symmetric regular Hadamard matrix of order $n^2$.
\end{theorem}

A list of $11$ classes of orders of symmetric Hadamard
matrices appears in \cite{Seberry:1972:LNM:292}, Appendix D. Some of these classes are infinite. For example class {\sc SHIII} is an infinite class, as a consequence of Dirichlet's theorem on the existence of primes in arithmetic progressions.

Another infinite series has been discovered recently in \cite[Lemma 1]{Balonin:Seberry:2015}:
\begin{theorem}
Let $q \equiv 1 \pmod{4}$ be a prime power. Then propus matrices exist for orders $n=\frac{q+1}{2}$ which give propus-Hadamard
matrices of order $2(q+1)$.
\end{theorem}

This series is derived from Turyn's infinite series of Williamson matrices \cite{Turyn:1972} and can be plugged into the GP array. Another construction of symmetric Hadamard 
matrices of the same order has been known for long time, 
see \cite[Lemma 5.2, p. 339]{Seberry:1972:LNM:292}.

One of us (D.D.) has subsequently observed that the same construction is applicable to the infinite series of Hadamard matrices of Goethals--Seidel type constructed by
Xia, Xia, Seberry, and Wu \cite[Theorem 3]{XXSW:2005}:
\begin{theorem}
Let $q=4n-1$ be a prime power $\equiv 3 \pmod{8}$. Then there  exists an Hadamard matrix of order $4n$ of Goethals--Seidel type in which
$$
(I-A)^T=-I+A, \quad B^T=B \quad \text{\it and} \quad C=D.
$$
\end{theorem}
In fact their matrix is not just a Hadamard matrix but also
a skew Hadamard matrix.

Instead of plugging the quadruple $[A,B,C,D]$ into the
Goethals--Seidel array, we can plug the permuted quadruple
$[B,C,D=C,A]$ into the GP array to obtain a propus-Hadamard
matrix. Thus we have the following theorem:
\begin{theorem}
Let $q=4n-1$ be a prime power $\equiv 3 \pmod{8}$. Then there  exists a symmetric Hadamard matrix of order $4n$ of GP-type,
i.e., obtained by using the GP array.
\end{theorem}

For $n=11,17,33,35,53,71,77,83,123,125$ these symmetric
Hadamard matrices are displayed on Balonin's webpages
\url{http://mathscinet.ru/catalogue/propus/dragomir/}.

\section{Overview of the algorithm for order 116}

We want to construct a SDS $[A,D,B,C]$ in $\bZ_{29}$ with 
parameters 
$(29;13,13,11,11;19)$ such that the subset $A$ is symmetric 
and $B=C$. This is necessary in order to use the GP array. 
Thus $A$ and $D$ have cardinality 13, and $B$ 
cardinality 11. 
(The other option, with $A$ and $D$ of cardinality 11 and
$B$ of cardinality 13, was treated separately in the same manner.)
We generate three files, one for each of $A,D,B$. The 
$A$-file contains the symmetric subsets of cardinality 13, 
the $D$-file arbitrary subsets of cardinality 13, and the 
$B$-file arbitrary subsets of cardinality 11. 
(Each subset is recorded on a separate line.) 
We do not collect all such subsets, but only those that 
pass the power spectral density (PSD) test. This test is an important tool as it 
cuts down the size of the file considerably. 
For a description of the PSD test see 
\cite[section 4]{DK:DCC:2014}.
For each of these three files and subset recorded there, we compute and record in a new file the periodic autocorrelation function (PAF) of the corresponding binary sequence of length 29. 
The subsets $A,D,B,C=B$ will form a SDS if and only if the 
sum of their PAFs takes the value $\lambda=19$ at all 
shifts different from 0. Thus the search for SDSs boils down 
to selecting one line in each of the $A,D,B$ PAF-files such 
that the element-wise sum of the first, second and twice the 
third line is equal to 19.

To find such a triple of lines, we adapted the meet-in-the-middle parallel collision finding technique of  \cite{vanOorschot:Wiener:JoC:1999}. Let us represent a line in the PAF file of $A$ as the sequence $(a_1 \enskip a_2 \cdots a_n)$, and similarly for the PAF files of $D$ and $B$. To find a triple of lines such that
\begin{equation}
a_i + d_i + 2b_i = \lambda,  \quad \forall \enskip i \in \{1, \ldots, n\},
\end{equation}
\noindent we define two functions, $f_{ad}$ and $f_b$:
\begin{eqnarray}
 f_{ad}(i, j) &:=& \text{Element-wise sum of line $i$ from $A$ and line $j$ from $D$}\\
 f_b(k) &:=& \text{Element-wise difference of $\lambda$ and twice line $k$ from $B$}
\end{eqnarray}

\noindent With these definitions, any case where  $f_{ad}(i, j) = f_b(k)$ constitutes a solution, and thus the existence of a symmetric Hadamard. 

To execute the search, we have a large number of processors perform random walks through the space of combinations of lines. Walks start at random positions, and deterministically decide at each step whether to execute $f_{ad}$ or $f_b$, and which lines of the file to read. To determine which function to perform in the next step, we concatenate the values in the result of $f_{ad}$ or $f_b$, and run that through a SHA-1 hash function. An indicator for the next function to perform, and corresponding line indices, are then derived from this hash value. A walk terminates when the resultant hash value reaches some pre-defined condition (usually a certain number of 0s at the beginning of the hash string, the choice of which depends on the time-memory tradeoff between computation and storage time). Starting and ending points of all completed walks are stored in a set shared by all processors.  This means that, with high probability, if two concatenated sums from both $f_{ad}$ and $f_b$ are the same, then the hash value will be the same, and the two walks will `merge' and arrive at the same point in the collective set of stored walks. When such an instance occurs, we have found our solution.

We implemented this algorithm in C++11, using Boost.MPI. The specific implementation was adapted from 
\cite{DiMatteo:2015}. All the data from the files was stored in a SQLite database. As the initial files are rather large, we perform a preprocessing step before doing the matching. All three initial PAF files are divided into subfiles based on the first number in each line. Then, only combinations of three subfiles such that the first numbers of each line sum to $\lambda$ are actually run through the program.

As the algorithm is random and parallel, it is difficult to benchmark its runtime. We ran it on SHARCNET's Orca cluster, a machine whose nodes have processors with speeds of either 2.2GHz or 2.7GHz, and a minimum of 32GB RAM. We used 16 MPI processes: 14 of them continuously executed random walks, one held the shared set of completed walks, and the last was responsible for receiving pairs of walks which ended at the same point, and extracting possible solutions. Searches were done for a fixed amount of time, usually 24h. Of the 15 matches found for $v = 29$, the shortest time taken was 238 seconds; the longest took over nine hours. All but five of the matches were found in less than three hours.

\section{Results}

In this section we present the construction of symmetric Hadamard matrices of orders
$4 \cdot 23 = 92$, $4 \cdot 29 = 116$ and $4 \cdot 43 = 172$.
The order $92$ is not new. The orders $116, 172$ are new.

The solutions are listed in the form of SDSs with four base blocks. From them one can construct the corresponding binary sequences and also the circulant matrices. To be specific, we label the positions of a binary sequence of length $v$ with $0,1,\ldots,v-1$ in that order. To a given subset
$X\subseteq\{0,1,\ldots,v-1\}$ we associate the binary sequence whose $-1$ entries occur exactly at the positions labeled by the elements of $X$. Further, to such a binary sequence we associate the circulant matrix of order $v$ whose first row is that sequence.
These circulant matrices can be plugged in to the GP array, in a suitable order, to obtain the symmetric Hadamard matrix. We say that a block of a SDS is
{\em symmetric} if the corresponding binary circulant matrix is symmetric.

\subsection{Four non-equivalent solutions for $v=23$}

\begin{figure}[h!]
 \centering
 \includegraphics[scale=0.4]{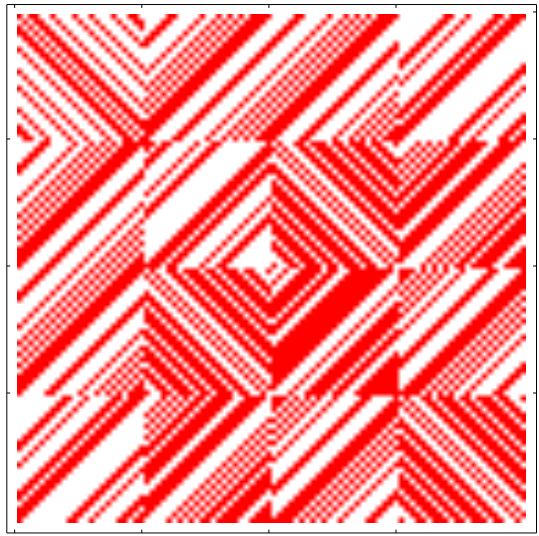}
 \caption{A visualization of the first solution of order 4$\cdot$23 = 92. White indicates a value of +1, and red -1. We can see clearly the symmetry of the matrix, as well as the distinct blocks of the GP array.}
 \label{hadamard41}
 \end{figure}

All four SDSs have parameters $(23;10,10,9,8;14)$ and are written as $[B,C=B,A,D]$ with $A$ symmetric. The quadruple of corresponding circulant matrices $[A,B,C,D]$ should be plugged in to the GP array. A graphical representation of the first solution listed below is displayed in Figure \ref{hadamard41}.

\begin{verbatim}
[[[0,1,2,3,5,7,9,12,17,18],[0,1,2,3,5,7,9,12,17,18],
     [0,2,3,6,10,13,17,20,21],[0,1,2,4,5,10,13,14]],
    [[0,1,2,3,5,7,10,11,13,19],[0,1,2,3,5,7,10,11,13,19],
     [0,5,7,8,11,12,15,16,18],[0,1,2,7,10,11,14,16]],
    [[0,1,2,3,6,8,9,10,14,19],[0,1,2,3,6,8,9,10,14,19],
     [0,1,3,8,11,12,15,20,22],[0,1,3,5,7,10,13,16]],
    [[0,1,2,5,6,8,10,13,14,16],[0,1,2,5,6,8,10,13,14,16],
     [0,1,3,4,10,13,19,20,22],[0,1,2,5,7,12,16,18]]]:
\end{verbatim}

\subsection{Fifteen non-equivalent solutions for $v=29$}

\begin{figure}[h!]
 \centering
 \includegraphics[scale=0.37]{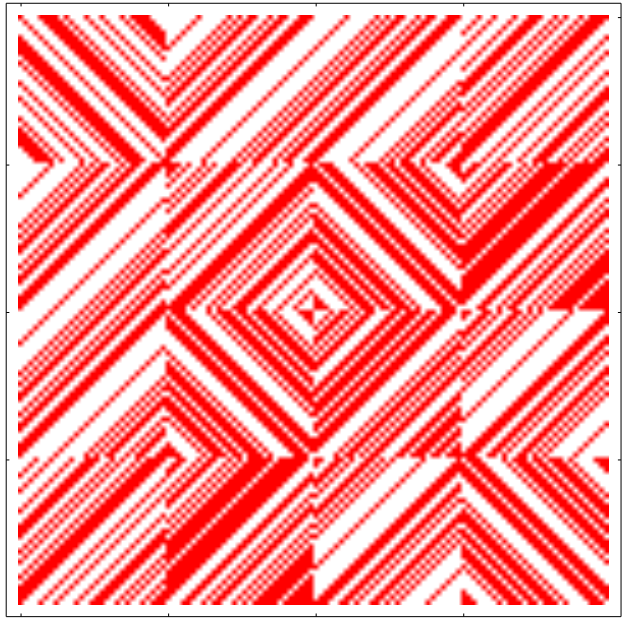}
 \caption{A visualization of the first solution of order 4$\cdot$29 = 116.}
 \label{hadamard42}
 \end{figure}

All 15 SDSs have parameters $(29;13,13,11,11;19)$. The first six are written as $[A,D,B,C=B]$ and the remaining nine as 
$[B,C=B,A,D]$, with $A$ symmetric in all cases. The quadruple of corresponding circulant matrices $[A,B,C,D]$ should be plugged in to the GP array. The first solution from those listed below is displayed in Figure \ref{hadamard42}.

\begin{verbatim}
[[[0,4,5,6,7,9,13,16,20,22,23,24,25],
     [0,1,2,3,5,8,10,13,14,15,18,22,25],
     [0,1,2,5,6,8,11,12,14,16,22],
     [0,1,2,5,6,8,11,12,14,16,22]],
    [[0,1,5,8,12,13,14,15,16,17,21,24,28],
     [0,1,4,5,6,7,10,11,14,16,17,19,24],
     [0,1,2,4,7,8,10,12,15,19,21],
     [0,1,2,4,7,8,10,12,15,19,21]],
    [[0,2,6,7,8,10,11,18,19,21,22,23,27],
     [0,1,3,4,5,7,10,13,14,15,17,20,24],
     [0,1,2,5,6,7,9,12,14,20,23],
     [0,1,2,5,6,7,9,12,14,20,23]],
    [[0,2,3,4,6,11,13,16,18,23,25,26,27],
     [0,1,2,5,6,8,10,11,14,17,18,20,22],
     [0,1,2,3,6,7,11,13,14,17,22],
     [0,1,2,3,6,7,11,13,14,17,22]],
    [[0,6,8,9,10,12,13,16,17,19,20,21,23],
     [0,1,3,5,6,7,10,12,13,16,18,21,22],
     [0,1,2,3,7,9,12,14,17,18,22],
     [0,1,2,3,7,9,12,14,17,18,22]],
    [[0,1,5,8,12,13,14,15,16,17,21,24,28],
     [0,1,2,3,6,8,9,12,13,15,19,20,23],
     [0,1,2,5,7,9,11,12,15,17,20],
     [0,1,2,5,7,9,11,12,15,17,20]],
    [[0,1,2,3,6,8,9,11,12,16,18,20,25],
     [0,1,2,3,6,8,9,11,12,16,18,20,25],
     [0,5,6,8,9,13,16,20,21,23,24],
     [0,1,2,3,4,7,8,11,13,17,19]],
    [[0,1,2,3,4,6,9,10,11,14,17,19,23],
     [0,1,2,3,4,6,9,10,11,14,17,19,23],
     [0,2,7,11,12,14,15,17,18,22,27],
     [0,1,2,3,6,9,13,14,18,20,24]],
    [[0,1,2,3,5,6,9,11,12,15,17,22,24],
     [0,1,2,3,5,6,9,11,12,15,17,22,24],
     [0,2,6,7,8,11,18,21,22,23,27],
     [0,1,2,3,5,6,10,13,14,17,21]],
    [[0,1,2,4,5,6,9,10,13,15,16,17,23],
     [0,1,2,4,5,6,9,10,13,15,16,17,23],
     [0,1,2,7,10,12,17,19,22,27,28],
     [0,1,3,5,7,10,13,16,19,21,25]],
    [[0,1,3,5,7,8,10,12,13,14,18,21,22],
     [0,1,3,5,7,8,10,12,13,14,18,21,22],
     [0,3,6,11,12,13,16,17,18,23,26],
     [0,1,2,3,4,6,10,12,16,17,20]],
    [[0,1,2,3,5,7,8,10,14,16,19,20,24],
     [0,1,2,3,5,7,8,10,14,16,19,20,24],
     [0,1,2,4,10,11,18,19,25,27,28],
     [0,1,4,5,7,9,12,13,16,20,23]],
    [[0,1,2,4,6,7,8,10,11,14,17,19,22],
     [0,1,2,4,6,7,8,10,11,14,17,19,22],
     [0,1,2,8,10,14,15,19,21,27,28],
     [0,1,2,5,7,10,14,15,18,19,24]],
    [[0,1,3,4,6,7,8,11,13,15,17,22,23],
     [0,1,3,4,6,7,8,11,13,15,17,22,23],
     [0,1,9,11,12,14,15,17,18,20,28],
     [0,1,2,5,6,10,12,17,18,21,26]],
    [[0,1,2,3,5,7,9,11,12,14,15,20,24],
     [0,1,2,3,5,7,9,11,12,14,15,20,24],
     [0,5,6,8,9,13,16,20,21,23,24],
     [0,1,4,5,7,8,10,16,17,18,23]]]:
\end{verbatim}

\subsection{One solution for $v=43$}

\begin{figure}[h!]
 \centering
 \includegraphics[scale=0.4]{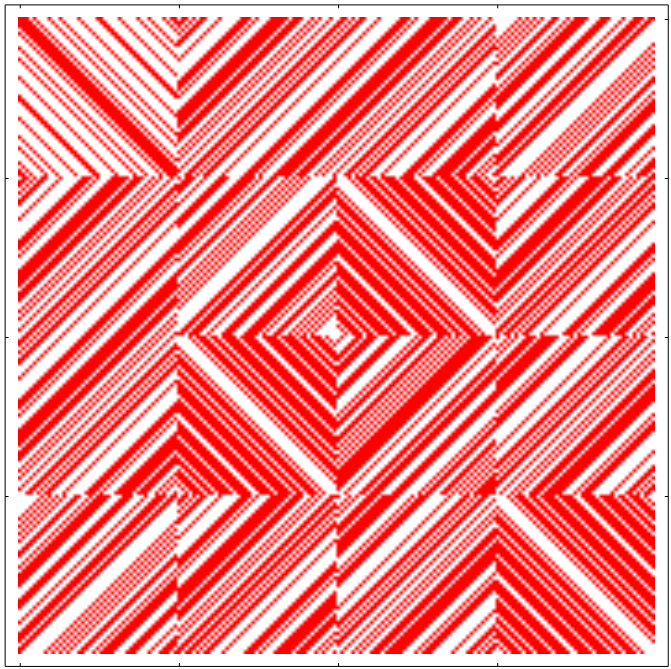}
 \caption{A visualization of the single solution found of order 4$\cdot$43 = 172.}
 \label{hadamard43}
 \end{figure}

In this case we start by constructing D-optimal matrices of order $86$ by using the well known two-circulant construction.
The parameter set of the relevant SDS, $[D,A]$, is $(43;21,15;15)$. The important feature of this SDS is that
one of the sets has to be symmetric. Our example is the
following:
\begin{verbatim}
[[0,1,3,4,5,8,12,13,14,18,19,20,21,23,26,27,29,30,32,34,36],
 [0,1,2,4,8,11,16,21,22,27,32,35,39,41,42]].
\end{verbatim}
The second block, A, is symmetric. When combined with the Paley difference set $B$ in $\bZ_{43}$, we can plug $[A,B,C=B,D]$ into the GP array (\ref{GP-array}) to obtain
a symmetric Hadamard matrix of order $4 \cdot 43 = 172$. This matrix is visualized in Figure \ref{hadamard43}.

\section{Acknowlegdements}
The authors wish to acknowledge generous support by NSERC.
This work was made possible by the facilities of the Shared
Hierarchical Academic Research Computing Network (SHARCNET) and Compute/Calcul Canada.
We thank N. A. Balonin and J. Seberry for sending us copies of their
papers \cite{Balonin:Seberry:2015,Balonin:Seberry:2014} and Michele Mosca for useful discussions.
O.D.M. would like to thank Gary Graham for some suggestions regarding
the adaptation of the parallel collision finding implementation.


\begin{thebibliography}{99}

\bibitem{Balonin:Seberry:2015}
N. A. Balonin, Jennifer Seberry, Visualizing Hadamard matrices:
the propus construction, preprint 15pp (submitted 6 Aug 2014).

\bibitem{Balonin:Seberry:2014}
N. A. Balonin, Jennifer Seberry, A review and new symmetric conference matrices, Informatsionno-upravliaiushchie sistemy, 2014, 8470; 4 (71), 2--7.

\bibitem{CraigenKharaghani:2007}
R. Craigen and H. Kharaghani, Hadamard Matrices and Hadamard
Designs. In Handbook of Combinatorial Designs. Edited by
Charles J. Colbourn and Jeffrey H. Dinitz. Second edition. Discrete Mathematics and its Applications (Boca Raton).
Chapman \& Hall/CRC, Boca Raton, FL, 2007.

\bibitem{DK:JCD:2012}
D. {\v{Z}}. {\Dbar}okovi{\'c} and I. S. Kotsireas,
New results on D-optimal matrices.
J. Combin. Designs, 20 (2012), 278--289.

\bibitem{DK:DCC:2014}
D. {\v{Z}}. {\Dbar}okovi{\'c} and I. S. Kotsireas, 
Compression of periodic complementary sequences and applications,
Des. Codes Cryptogr. 74 (2015), 365--377.

\bibitem{IoninShrikhande:CUP:2006}
Y. J. Ionin and M. S. Shrikhande, Combinatorics of Symmetric
Designs. New Mathematical Monographs, 5. Cambridge University Press, Cambridge, 2006.

\bibitem{Mathon:1978}
R. Mathon, Symmetric conference matrices of order $pq^2 +1$.
Canad. J. Math., 30 (1978), 321--331.

\bibitem{Seberry:1972:LNM:292}
J. Seberry Wallis, Hadamard Matrices,
in W. D. Wallis, A. Penfold Street, Jennifer Seberry Wallis, Combinatorics: Room squares, sum-free sets, Hadamard matrices. Lecture Notes in Mathematics, Vol. 292.
Springer-Verlag, Berlin-New York, 1972.

\bibitem{Turyn:1972}
R. J. Turyn,
An infinite class of Williamson matrices,
J. Combinatorial Theory Ser. A 12 (1972), 319--321.

\bibitem{vanOorschot:Wiener:JoC:1999}
Paul C. van Oorschot and Michael J. Wiener,
Parallel collision search with cryptanalytic applications,
Journal of Cryptology, January 1999, Volume 12, Issue 1, 1--28.


\bibitem{XXSW:2005}
M. Xia, T. Xia, J. Seberry and J. Wu,
An infinite series of Goethals--Seidel arrays,
Discrete Applied Mathematics 145 (2005) , 498--504.

\bibitem{DiMatteo:2015}
O. Di Matteo, \emph{Parallelizing quantum circuit synthesis.} MSc thesis,
University of Waterloo (2015).
\end{thebibliography}
\end{document}